\documentclass[11pt]{amsart}

\usepackage{graphicx}
\usepackage{amsfonts, amsmath, amsfonts, amssymb}
\newtheorem{thm}{Theorem}[section]
\newtheorem{cor}[thm]{Corollary}

\newtheorem{prop}[thm]{Proposition}
\theoremstyle{definition}
\newtheorem{defn}[thm]{Definition}
\theoremstyle{remark}
\newtheorem{rem}[thm]{Remark}
\numberwithin{equation}{section}
\theoremstyle{quest}

\numberwithin{equation}{section}
\theoremstyle{prob}

\numberwithin{equation}{section}
\theoremstyle{answer}

\numberwithin{equation}{section}
\theoremstyle{conj}
\newtheorem{conj}[]{Conjecture}
\numberwithin{equation}{section}

\begin{document}

\title[PLembedd-and-qr]{On a construction of Burago and Zalgaller}
\author{Emil Saucan}%
\address{Department of Mathematics, Technion, Haifa, Israel}%
\email{semil@tx.technion.ac.il}%

\thanks{Research supported by the Israel Science Foundation Grant 666/06 and by European Research Council under the European Community's Seventh Framework Programme
(FP7/2007-2013) / ERC grant agreement n${\rm ^o}$ [203134].}%

\date{\today}

\begin{abstract}
The purpose of this note is to scrutinize 
the proof of Burago and Zalgaller regarding the existence of $PL$ isometric embeddings of $PL$ compact surfaces into $\mathbb{R}^3$. We conclude that their proof does not admit a direct extension to higher dimensions. Moreover, we show that, in general,  $PL$ manifolds of dimension $n \geq 3$ admit no nontrivial 
$PL$ embeddings in $\mathbb{R}^{n+1}$ that are close to conformality. We also extend the result of Burago and Zalgaller to a large class of noncompact $PL$ 2-manifolds.
The relation between intrinsic and extrinsic curvatures is also examined, and we propose a $PL$ version of the 
 Gauss compatibility equation for smooth surfaces. 

\end{abstract}
\maketitle


\section{Introduction and Main Results}

In \cite{bz} Burago and Zalgaller proved the following theorem, that represents 
a $PL$ version for dimension $n = 2$ of the celebrated Nash-Kuiper $\mathcal{C}^1$ isometric embedding theorem \cite{na1}, \cite{ku}:

\begin{thm} \label{thm:BZ}
Any compact orientable $PL$ $2$-manifold admits a $PL$ isometric embedding in $\mathbb{R}^3$.
\end{thm}

\begin{rem}
Nonorientable $PL$ 2-manifolds are shown to admit $PL$ immersions into $\mathbb{R}^3$.
\end{rem}

Of course, one has to 
properly define the notion of isometric embedding for the case of
$PL$ manifolds. We leave this for Section \ref{sec:Terminology and notation}.

The main purpose of this note is to examine the validity of Theorem \ref{thm:BZ} above in
dimensions $n > 2$, hence 
of the Nash-Kuiper-Burago-Zalgaller embedding process 
 - henceforward abbreviated as $NKBZ$. In particular 
we prove the following negative result:

\begin{thm} \label{thm:main}
In any dimension $n \geq 3$ there exists a compact $PL$ manifold (in fact an infinity of such manifolds) that can
not be $PL$ isometrically embedded in $\mathbb{R}^{n+1}$ via 
the NKBZ method.
\end{thm}

The main ingredient in the proof of this
theorem, besides a scrutiny of the proof of Theorem  \ref{thm:BZ}, consists in computing the {\it coefficients of conformality} (see Section 2.2 below) of a certain standard mapping of the neighborhoods of the vertices, that appears in the Burago-Zalgaller construction (see Section 3).

We can, in fact, strengthen the result above, 
as follows:

\begin{prop}
In any dimension $n \geq 3$ there exists a compact $PL$ manifold (in fact an infinity of such manifolds) that admits 
no nontrivial 
quasiconformal 
embedding in $R^{n+1}$. 
\end{prop}

Again, as in the proof of Theorem 1.3, the method to prove the result above,
besides an examination of the proof of Theorem  \ref{thm:BZ}, is to make appeal to the theory of quasiconformal/quasiregular mappings, more precisely to apply a theorem of Rickman and Srebro (\cite{rs}) on the nonexistence of quasiregular mappings with large local index on a finite, evenly distributed set in $\mathbb{R}^{n}, n \geq 3$ (see Theorem \ref{thm:RS} below).


\begin{rem}
The implication of the nonexistence of such an embedding to Graphics, Imaging and related applicative fields was discussed in some detail in \cite{Sa10}.
\end{rem}

The remainder of the paper is apportioned as follows: In Section 2 we present the necessary background. Section 3 represents a sketch of the Burago-Zalgaller construction. It is followed, in Section 4, by an analysis of the obstructions to the extension of the said 
construction to higher dimensions. The main results are proven in Section 5. 
In the last section we discuss the role of curvature. In particular, we show that the compactness condition, given in \cite{bz}, is too restrictive,  and prove that their embedding result holds, in fact, for a quite large class of unbounded manifolds (with or without boundary).
We also give a $PL$ version 
of the Gauss compatibility condition for existence of isometric embeddings of $PL$ 2-manifolds into $\mathbb{R}^3$.

A few precautionary words to 
the reader: It is possible -- indeed, it is rather probable -- 
that part of the material contained herein, especially in Section 2, will appear redundant and classical. 
However, since the paper does not properly appertain neither to Differential Geometry, nor to Quasiconformal/Quasiregular
Function Theory, but rather lies in an indeterminate area between these two fields, we have decided to make the paper self contained, 
hence as friendly as possible to readers of both backgrounds (and hopefully, of any
mathematical background).


\section{Terminology, Notation and Preliminaries} \label{sec:Terminology and notation}

\subsection{$PL$ Isometric Embeddings}

It is a quite common mistake to believe $PL$ isometric embeddings coincide with the isometric embeddings in the classical (smooth) Riemannian context.
It is true that a (rather straightforward) Riemannian structure on $PL$ manifolds can be defined 
-- for full details see \cite{Te}.
However, due to the lack of smoothness, the classical curvature operator can not be defined, therefore the two notions of isometric embedding (i.e. classical and $PL$) diverge. Indeed, they coincide only for piecewise flat manifolds -- see, e.g. \cite{Be}. It is, therefore important to emphasize the difference and give 
the correct definition in the $PL$ case. (See \cite{gr2}, \cite{LeDo} and, for a a lengthier discussion, \cite{Sa10}).

\begin{defn}
Let $M^n$ be a $PL$ manifold (or, more generally, a metric space). A map $f:M^n \rightarrow N^\nu$, where $N^\nu$ is another manifold (metric space) is called a $PL$ {\it isometric embedding}, or a {\it path isometry} iff

(i) It is a topological embedding

and

(ii) \({\rm length}(f(\gamma)) = {\rm length}(\gamma)\,,\)

for all\footnote{even for those of infinite length} 
curves in $M^n$.

\end{defn}

\begin{rem}
It is interesting to note that if one discards the embedding requirement, then the following rather surprisingly
result is obtained, namely that any $PL$-manifold of dimension $\leq 4$ admits an $PL$ isometry into
$\mathbb{R}^n$ (\cite{za}).

For more results on the metric geometry of $PL$ surfaces see e.g. \cite{Sa04}, \cite{sa}.
\end{rem}


\subsection{Locally Flat Manifolds}

In the following 
 we shall want to clearly distinguish between piecewise flat embeddings and simply 
$PL$ ones. Therefore, we bring here, following Loukkainnen \cite{Lu}, the necessary definition.
First, we have to introduce some notation:

Let $m \neq n$ be natural numbers. First, let $\mathbb{R}^n_+ = \{X \in \mathbb{R}^n\,|\, x_n \geq 0\}$, $\mathbb{R}^n_{++} = \{x \in \mathbb{R}^n_+\,|\, x_{n-1} \geq 0\}$, and we identify $\mathbb{R}^q$ with $\{x \in \mathbb{R}^n\,|\, x_i = 0, {\rm if}\; i \geq q+1\}$. Next, we define: $\mathbb{R}^{n,q}_+ = \{x \in \mathbb{R}^n_+\,|\, x_i = 0, {\rm if}\; i \leq n-q\}$, and $\mathbb{R}^{n,q}_{++} = \{x \in \mathbb{R}^{n,q}_+\,|\, x_{n-1} \geq 0\}$. ($\mathbb{R}^{n,1}_{++} = \mathbb{R}^{n,1}_+$.)

A model for locally flat manifold pair at a point of a submanifold, or simply a $(n,q)$-{\it model}, we mean one of the following pairs:

(i) $(\mathbb{R}^n, \mathbb{R}^q)$, $(\mathbb{R}^n, \mathbb{R}^q_+)$, $(\mathbb{R}^n_+,\mathbb{R}^{n,q}_+)$, $(\mathbb{R}^n_+,\mathbb{R}^{n,q}_{++})$, if $q \geq 2$;

(ii) $(\mathbb{R}^n, \mathbb{R}^1)$, $(\mathbb{R}^n, \mathbb{R}^1_+)$, $(\mathbb{R}^n_+,\mathbb{R}^{n,1}_+)$, if $q = 1$;

(iii) $(\mathbb{R}^n, \mathbb{R}^0)$, $(\mathbb{R}^n_+, \mathbb{R}^0)$, if $q = 0$;

\begin{defn}
A $PL$ submanifold $M$ of $\mathbb{R}^n$ is called {\it locally flat} ($LF$), or {\it piecewise flat} ($PF$), iff for any point $x \in M$, there exists an open neighbourhood $U$ of $x$ in $\mathbb{R}^n$, an $(n,q)$-model $(N,L)$ and a $PL$ homeomorphism $h:U \rightarrow N$, such that $h(U \cap M) = L$, and $h(x) = 0$.
\end{defn}


\subsection{Quasiregular Mappings}

\begin{defn}
Let $D \subseteq \mathbb{R}^n$ be a domain; $n \geq 2$ and let $f: D \rightarrow \mathbb{R}^n$ be a continuous
mapping. $f$ is called
\begin{enumerate}
\item {\it quasiregular} (\it qr) iff
\begin{enumerate}
\item
is locally Lipschitz (and thus differentiable a.e.); 
\\ \hspace*{-0.8cm} and
\item \(0 < |f'(x)|^n \leq KJ_f(x)\), for any \(x \in M^n\);
\end{enumerate}
\;where $f'(x)$ denotes the formal derivative of $f$ at $x$, $|f'(x)| = \sup
\raisebox{-0.25cm}{\mbox{\hspace{-0.75cm}\tiny$|h|=1$}}|f'(x)h|$, and where $J_f(x) = detf'(x)$;
\item {\it quasiconformal} (qc) iff $f:D \rightarrow f(D)$ is a quasiregular homeomorphism;
\end{enumerate}
The smallest number $K$ that satisfies condition (b) above is called the {\it outer dilatation} of
\nolinebreak[4]$f$.
\end{defn}

\begin{rem}
One can extend the definitions above to mappings $f:M^n \rightarrow N^n$, where $M^n,N^n$ are oriented, connected
 Riemannian $n$-manifolds, $n \geq 2$, by using coordinate charts (for details see, e.g. \cite{v}).
\end{rem}

\begin{rem} \label{rem:bound-dist}
Admittedly, the appellative 
{\it quasiconformal} conveys little geometric meaning. Certainly, the term {\it mappings of bounded distortion}, originally introduced by Reshetnyak \cite{Re} is far more apt
to convey the geometric content 
of such mappings and would, perhaps, be better used in 
the present geometric context. However, the notions of {\it quasiconformal} and {\it quasiregular} mappings have gained ground, due to the fact that they clearly point 
to the fact that such mappings represent natural extensions of the classical conformal, respective analytic mappings in the plane. Because of this reason, and due to the fact that core of the paper is based upon the observation that no ``almost conformal'', quasiconformal parametrizations of the neighborhoods of the vertices of a polyhedron are possible in dimension $n \geq 3$ 
(see Section 4,\,(5) and Section 5), we have retained here the more established terminology (albeit, perhaps, at the detriment of geometric intuition).
\end{rem}

If $f: D \rightarrow \mathbb{R}^n$ is quasiregular, then there exists $K' \geq 1$ such that the following
inequality holds a.e. in $M^n$:

\begin{equation} \label{eq:qr1}
J_f(x) \leq K'\inf_{|h|=1}{|T_xfh|^n}
\end{equation}

By analogy with the outer dilatation we have the following definition:

\begin{defn}
The smallest number $K'$ that satisfies inequality (\ref{eq:qr1}) is the {\it inner dilation} $K_I(f)$ of $f$, and $K(f)
= \max(K_O(f),K_I(f))$ is the {\it maximal dilatation} of $f$. 
If $K(f) < \infty$ we say that $f$ is called $K$-qr.
\end{defn}

The dilations are $K(f), K_O(f)$ and $K_I(f)$ are simultaneously finite or infinite. Indeed, the following
inequalities hold: $K_I(f) \leq K_O^{n-1}(f)$ and $K_O(f) \leq K_I^{n-1}(f)$. This allows us to talk about the dilatation ($K$) of a mapping, without being more specific.

For planar quasiregular mappings we have the following classical structure theorem:

\begin{thm}[Stoilow, \cite{St}]
Let $f:D \rightarrow \mathbb{R}^2$ be a nonconstant quasiregular mapping. Then $f$ admits the following factorization: $f = g \circ h$, where $h:D \rightarrow \mathbb{R}^2$ is quasiconformal and $h:f(D) \rightarrow \mathbb{C} \equiv \mathbb{R}^2$ is a nonconstant holomorphic 
mapping.
\end{thm}

In particular, any quasiregular mapping $f:D \rightarrow \mathbb{R}^n$ is locally quasiconformally equivalent to the mapping $z \mapsto z^m, z \in \mathbb{C}$, for some $m = m(x_0) \in \mathbb{N}, x_0 \in D$.
Unfortunately, no generalization of Stoilow's Theorem result exists in dimension $n \geq 3$. (We shall further expand upon this subject shortly.)

\begin{defn}
Let $f:D \rightarrow \mathbb{R}^n$ be a quasiregular mapping. The set $B_f = \{x \in M^n\,|\, f\; {\rm is\; not\;
a\;
local\; homeomorphism\; at\;} x\}$ is called the branch set 
of $f$.
\end{defn}

Since any quasiregular mapping $f:M^n \rightarrow N^n$ is {\it discrete}, that is $f^{-1}(y)$ is discrete, for any $y \in M^n$ (see \cite{Re}), we can introduce the following

\begin{defn}
Let $f:D \rightarrow \mathbb{R}^n$ be an orientation preserving map. The {\it local topological index} of $f$ at $x$ is
defined as:
\begin{equation}
i(x,f) = \inf\raisebox{-0.25cm}{\hspace{-0.9cm}\mbox{\tiny $U \in \mathcal{N}(x)$}}\sup\raisebox{-0.25cm}{\hspace{-0.4cm}\mbox{\tiny ${y}$}}\;\,|f^{-1}(y) \cap U|
\end{equation}
\end{defn}

Note that if $f:D \rightarrow \mathbb{R}^n$ is a quasiregular mapping, then $i(x,f) \geq 1$ and, moreover, $x \in
B_f$ iff $i(x,f) > 1$.

Also, for $n \geq 3$ 
the local topological index cannot be uniformly too large on all the points of a non-degenerate continuum $F$. To
be more precise, the following inequality holds:
\begin{equation}
\inf\raisebox{-0.25cm}{\hspace{-0.65cm}\mbox{\tiny $x \in F$}}i(x,f) < n^{n-1}K_I(f)\,.
\end{equation}
 (See e.g. \cite{ric2}, III. 5.9.)

Moreover, even though local topological index can be arbitrarily large at an isolated point (see \cite{ric1}, pp.
263-264), it can not be too large even on a finite number of points if the points and the indices of the map $f$ at these points are evenly distributed: 

\begin{thm} [\cite{rs}, Theorem 1.1] \label{thm:RS}
Let $f:D \rightarrow \mathbb{R}^n,\, D \subseteq \mathbb{R}^n,\, n \geq 3$ be a non-constant $K$-$qr$ mapping.
Then, for any $x_0 \in D$, there exist $t_0, p > 0$, such that, and for any $x_1,\ldots, x_m \in
\mathbb{B}^n[x_0,t]$, $0 < 0 \leq t_0$\,, and satisfying the following conditions:
\begin{enumerate}
\item $|x_0 - x_m| = t$;
\item $|x_{j-1} - x_j| = t/p$, where $p_0 \leq p \leq m \leq p^\nu$ and  $1 \leq \nu \leq \big(i(x_0,f)/K_I(f)\big)^{\frac{1}{n-1}}$\,,
\end{enumerate}
there exists $j \in \{1,\ldots, m\}$ such that $i(x_j,f) < i(x_0,f)$. (Here $\mathbb{B}^n[x_0,t]$ denotes the
closed ball of center $x_0$ and radius $t$.)
\end{thm}

\subsubsection{Coefficients of Quasiconformality} \label{sec:coefCQ}

We bring below a few results regarding the coefficients of quasiconformality of some specific domains. In this we
rely on \cite{v} and \cite{car}. These results are needed in the proof of our result.

\begin{defn}
Let $D \subset \mathbb{R}^n, D \simeq \mathbb{B}^n$. The {\em coefficients of quasiconformality} of $D$ are
defined as follows:
\begin{equation}
K_I(D) = \inf_{f:D \stackrel{\!\!\sim}{\rightarrow}\mathbb{B}^n}K_I(f), \:\:\: K_O(D) = \inf_{f:D
\stackrel{\!\!\sim}{\rightarrow}\mathbb{B}^n}K_O(f),
\end{equation}

\[ K(D) = \inf_{f:D
\stackrel{\!\!\sim}{\rightarrow}\mathbb{B}^n}K(f).\]
\end{defn}

\begin{defn}
Let $x \in \mathbb{R}^n$ be a point with cylindrical coordinates $x = (r\cos\varphi,r\sin\varphi,z_1,\ldots ,
z_{n-2})$. The set $D_\alpha = \{0 < \varphi < \alpha\}$, $(0 < \alpha \leq 2\pi)$ is called a {\em wedge} of
angle $\alpha$. More generally, a domain $D \subset \mathbb{R}^n, n \geq 3$ is called a {\em wedge} of angle
$\alpha$ iff there exists a similarity mapping $f$ such that $f(D) = D_\alpha$. $f^{-1}(\{r=0\})$ is called the
{\em edge} of $D$. Given a domain $\Omega$, a point $b \in \partial \Omega$ is called a {\em wedge point} iff
there exist a neighborhood $U$ of $b$ and a wedge $D$ of angle $\alpha$, such that $b$ lies on the edge of
$\Omega$ and $U \cap D = U \cap \Omega$.
\end{defn}

\begin{defn}
The homeomorphism $f:D_\alpha \rightarrow D_\beta$, $f(r,\varphi,z) = (r,\frac{\alpha}{\beta}\varphi,z)$, $z = (z_1,\ldots ,z_{n-2})$, is called
a ({\em folding}, or {\em winding} (or, more precisely, ($k$){\em-winding} mapping).
\end{defn}

We should note that, for $n \geq 3$, foldings are topologically equivalent to $z \mapsto z^k \times Id$ (where $Id$ denotes the identity mapping of $\mathbb{R}^{n-2}$).

If $\alpha \leq \beta$, then $f$ is quasiconformal, with dilatations $K_I(f) = \frac{\alpha}{\beta}, K_O(f) \geq
(\frac{\alpha}{\beta})^{1/(n-1)}$. In particular, for $\beta = \pi$, we obtain $K_I(D_\alpha) =
\frac{\pi}{\alpha}, K_O(D_\alpha) = (\frac{\pi}{\alpha})^{1/(n-1)}$, whence $K(D_\alpha) = \frac{\pi}{\alpha}$.

\begin{rem}
Remarkably, the coefficients of quasiconformality for non-convex domains (i.e. $\pi \leq 2\pi$) are not known (at least to the best of our knowledge).
\end{rem}

Following \cite{car}, we note the following natural generalization of the definition of a wedge:

\begin{defn}
The domain $D_{\alpha k} \subset \mathbb{R}^n, D_{\alpha k} = \{(r,\varphi_1,\dots,\varphi_{n-k
+1},z_{n-k+1},\dots,z_n)\}$ is called a {\em dihedral wedge of type $k$ and angle $\alpha$}.
\end{defn}

\begin{rem}
For $k = n - 2$ we recuperate the classical definition of wedges.
\end{rem}

\begin{prop}[\cite{car}]
The coefficients of quasiconformality for $D_{\alpha k}$ are:
%
%
%
%
\begin{equation}
K_I(D_{\alpha k}) = \frac{\pi^{n-k-1}}{\alpha_1\cdots \alpha_{n-k-1}}\,\,,
\:\:
K_O(D_{\alpha k}) \geq \left(\frac{\pi^{n-k-1}}{\alpha_1\cdots \alpha_{n-k-1}}\right)^{\frac{1}{n-1}}\,,
\end{equation}

\[K(D_{\alpha k}) = \frac{\pi^{n-k-1}}{\alpha_1\cdots \alpha_{n-k-1}}\,\,.\]

\end{prop}

\begin{cor}
Let $D$ be a convex polyhedron in $\mathbb{R}^n$ and let $m$ denote the number of faces of $D$. Then we have the following estimates:
\begin{equation}
K_I(\mathcal{P}) \geq \frac{m-n+2}{m-n}\,\,,\:\: K_O(\mathcal{P}) \geq \left(\frac{m-n+2}{m-n}\right)^\frac{1}{n-1},
\end{equation}

\[K(\mathcal{P}) \geq \frac{m-n+2}{m-n}\,\,.\]

\end{cor}

\begin{rem}
Clearly, the same estimates hold for $PL$-smooth convex manifolds.
\end{rem}
A different slight generalization of wedges is the following one:

\begin{defn}
Let $D \subset \mathbb{R}^n$ be a domain. We say that $D$ has a {\em curvilinear wedge} of angle $\alpha$ at $x_0
\in \partial D$ iff, for all $K > 1$, there exists a neighbourhood $U$ of $x_0$, such that $f(U \cap D) =
\mathbb{B}^n \cap D_\alpha$.
\end{defn}

Another type of related domains are the so called {\em raylike} domains:

\begin{defn}
A domain $D \subset \mathbb{R}^n$ is called {\em raylike, with vertex} $v \in \partial D$, iff $v + t(x - v) \in
D$, for all $x \in D$ and any $t > 0)$.
\end{defn}

\begin{thm}[\cite{v}, Theorem 40.3] \label{thm:ray}
Let $D, G$ be domains in $\mathbb{R}^n$, such that $G$ is raylike, with vertex $v$. If $v$ has a neighbourhood
$U$, such that $U \cap D = U \cap G$, then $K_I(D) \geq K_I(G), K_O(D) \geq K_O(G), K(D) \geq K(G)$.


\end{thm}

\begin{rem}
These (rather straightforward) generalizations of the notion of wedge, allow us to extend Theorem 1.2 and Proposition 1.4 to more general $PL$ embeddings, 
not just to piecewise flat 
ones (see also Section 4,\,(3) below).
\end{rem}

Before concluding his section, let us consider again the folding map, this time from a different point of view:
As we have already noted, Stoilow's Theorem does not hold in dimensions higher than 2. There exists, however, a characterization of those quasiregular mappings topologically equivalent to a folding:

\begin{thm}[Martio-Rickman-V\"{a}isal\"{a}, \cite{MRV}]
Let $f:D \rightarrow \mathbb{R}^n$ be a quasiregular mapping and let $x_0$ be a point of $D$.
If there exist a neighbourhood $N_{x_0}$ of $x_0$, and a homeomorphism $\varphi_{x_0}:N_{x_0} \rightarrow \mathbb{R}^n$, such that $\varphi(N_{x_0} \cap B_f) \subset  \mathbb{R}^{n-2} \subseteq \mathbb{R}^n$,
then locally at $x_0$, $f$ is quasiconformally equivalent  to a folding mapping (hence to $z \mapsto z^k \times Id$).
\end{thm}



\begin{rem}
This result suggests an 
approach to the existence problem of $PL$ embeddings of $PL$ $n$-manifolds into $R^{n+1}$, $n \geq 3$, alternative to the one adopted in the present paper.
\end{rem}



\subsection{Higher Dimensional Angles} 

Since we discuss in the following sections -- mainly in Sections 5 and 6 -- the role of angles (and curvatures) in the case of dimension higher than 2, we succinctly present here a modicum of necessary definitions.

While presumably intuitive, we bring here the following technical definition of dihedral angles, as given in \cite{cms}:

\begin{defn}
A {\em simplicial con}e $C^k \subset \mathbb{R}^k \subset  \mathbb{R}^n$, is defined as: $C^k =
\bigcap_{\raisebox{-0.25cm}{\hspace{-0.5cm}\mbox{\scriptsize${j=1}$}}}^{\raisebox{0.25cm}{\hspace{-0.3cm}\mbox{\scriptsize${k}$}}}
\!\!H_j$\,, where $H_j$ are open half spaces in {\em general position}, such that $0 \in H_j\,, j = 1,\ldots,k$.
$L^{k-1} = C^k \bigcap \mathbb{S}^{n-1}$ is called a {\em spherical simplex}.
\end{defn}

\begin{defn}
Consider the simplices $\sigma^k < \tau^m$, and let $p \in \sigma^k$. Define the {\em normal cone}:
$C^{\bot}(\sigma^k,\tau^m) = \{\overrightarrow{px}\,|\, x \in \tau^m, \; \overrightarrow{px} \bot\, \sigma^k\}$,
where  $\overrightarrow{px}$ denotes the ray through $x$ and base-point $p$. The spherical simplex
$L(\sigma^k,\tau^m)$ associated to $C^{\bot}(\sigma^k,\tau^m)$  is called the {\em link} of $\sigma^k$ in
$\tau^m$.
\end{defn}

\begin{rem}
$C^{\bot}(\sigma^k,\tau^m)$ does not depend upon the choice of $q$.
\end{rem}

\begin{defn}
The {\em (internal) dihedral angle} $\measuredangle(\tau^k,\sigma^m)$ is defined as the normalized volume of
$L(\sigma^k,\tau^m)$, where the normalization is such that the volume of $\mathbb{S}^{n-1}$ equals $1$, for any $n
\geq 2$.
\end{defn}

\begin{defn}
Denote by $L(\sigma^k,\tau^m)*$ the {\em dual simplex} of $L(\sigma^k,\tau^m)$, i.e. $L(\sigma^k,\tau^m)* = \{v
\in \mathbb{S}^{} \,|\, \measuredangle(v,u) > \pi/2, \forall u \in  L(\sigma^k,\tau^m) \}$. The normalized volume
of $L(\sigma^k,\tau^m)*$ is called the {\em exterior dihedral angle} of $\sigma^k \subset \tau^m$ and we denote it
by $\angle(\tau^k,\sigma^m)$.
\end{defn}



\section{The Burago-Zalgaller construction}

We present here very briefly the main steps of the Burago-Zalgaller construction, the accent being placed on the geometric aspects of the construction and on those elements of the proof that are problematic when passing to higher dimensions -- to be discussed in detail in the following section. For the full technical 
intricacies of the proof, the reader should consult the original paper \cite{bz}. We concentrate solely on the case of compact, orientable manifolds, both because they represent the basic case (whose modification produces the construction for the other cases) and because the full connection with quasiconformal mappings is displayed here (in contrast with the nonorientable case).

\begin{enumerate}

\item Start with a $\mathcal{C}^2$ smooth, short embedding $f_0$ of the given 
$PL$ compact, closed 2-manifold $\mathcal{P}$. To obtain the necessary short embedding, one may use a suitable 
    homotety (see \cite{na1}).

\item Produce a variation of the original embedding in certain (small) neighborhoods of the vertices, such that each neighborhood has a standard form (more precisely, a disk neighborhood) that allows us to produce a {\it standard embedding} in the vicinity of the vertices.
    The standard embedding above is supplied by the {\it standard conformal map} (or {\it folding}) from $K(\theta,\rho) =
\{0 \leq \varphi \leq \theta, \rho > 0\}$ to $K(\lambda,r) = \{0 \leq \psi \leq \lambda, r >
0\}$ given by:

\[\psi = \frac{\lambda}{\theta}\varphi,\; r = a\rho^{\lambda/\theta}.\]

(The most important case for our purposes being: $\lambda = 2\pi$.)

    The resulting, varied embedding $f_1$ will have different forms if the sum $\theta_i$ of the angles incident  at the vertex $v_i$ is $< 2\pi$ or $> 2\pi$ -- see also (4) below. 

   However, in both cases, the embedding $f_1$ will be short in the complement of the said neighborhoods.

\item Replace the neighborhoods of the vertices with small disjoint polygonal neighborhoods $N_v$. More precisely:
  \begin{enumerate}
           \item If $\theta < 2\pi$, encircle $A$ by a small ``regular'' hexagon composed of $6$ triangles of apex angle $\theta/6$.

          Some small enough neighbourhood of $A$ the will be mapped by the standard conformal mapping onto a planar disk.

          Over each triangle included in such a neighbourhood, one can perform the basic construction (see (6) below), obtaining a $PL$ isometric embedding of this neighbourhood.
           \item If $\theta > 2\pi$, proceed analogously to the previous case but

           \begin{enumerate}
           \item In a small circular neighbourhood of radius $r_1$ map (a) isometrically on radial
segments and (b) using a $\theta/2\pi$ contraction on circles centered at $A$;
           \item In a annular neighbourhood $\{r_1 < r < r_2\}$ use the standard conformal mapping with
the same contraction factor $\theta/2\pi$.
           \end{enumerate}
     Replace the neigbourhood above with a ``cogwheel'' (i.e. a circle surrounded by isosceles ``triangles'' of sides, e.g. $2\delta$, and having as bases arcs of the same length).
     %
     %
     The interior of each ``cogwheel'' is $PL$ isometric embedded using ``ripples''. (The basic element of each such ``ripple'' is a pair of congruent triangles, having a common vertex in the center of the ``cogwheel'', one side of each being a radius, and a second common vertex built over the midpoint of an arc used in the construction of the ``cogwheel'' -- see Figure 4 of \cite{bz}).
  \end{enumerate}

The complement of 
    $\bigcup_vN_v$
    is triangulated using solely {\it acute angle triangles}. (In particular, at convex vertices subdivide each triangle into $n^2$ similar triangles, for some large enough $n$; while at non-convex vertices into almost regular triangles.) This has to be performed with care, so that a certain inequality (\cite{bz}, (5), p. 373), regarding the angles of the triangles composing the triangulation would(will) hold. Moreover, an additional variation of the triangulation is also applied.

\item A further variation of $f_1$ is performed in the complement of the union of the neighborhoods of the vertices.
    Here is employed a construction of Kuiper \cite{ku} constituting in the adding a (finite) succession of $\mathcal{C}^2$ smooth waves superimposed on $f_1(\mathcal{P})$. 
    The embedding $f_2$ obtained in this manner is $\mathcal{C}^2$ smooth and short\footnote{Again, the contraction constant can be precisely specified.} Moreover, $f_2$ will be {\it almost conformal} in the complement of (sufficiently small) neighborhoods of the vertices $v_i$ such that the $\theta_i > 2\pi$. (The measure of ``almost conformality'' is a function of $\alpha$ and the degree of approximation of an isometry.)

\item The triangles obtained at the previous stages are further subdivided, in order to obtain of sufficiently small mesh. Here the degree of ``almost conformality'' is exceedingly important, since it is used to assure that the triangles $T_i$ and $t_i$ (of the {\it standard construction element} -- see (6) below) are {\it almost similar} and also to produce dihedral angles close to $\pi$. In addition, the dihedral angle between adjacent ``$t_i$'' triangles 
    are uniformly close to $\pi$.
\item Apply the {\it canonical} ({\it standard}) {\it construction element}:

      \begin{enumerate}
      \item Let $T = \triangle(A_1,A_2,A_3)$ and $t =
 \triangle(a_1,a_2,a_3)$
be acute triangles;
      \item let $B,b$ and $R,r$ the centers and radii of their
respective circumscribed circles;
      \item let $E_p = \frac{1}{2}A_kA_l, e_p = \frac{1}{2}a_ka_l;\, p,k,l \in
\{1,2,3\}$;
      \item and let $H_p = BE_p, h_p = be_p$.

Moreover, let $T \simeq t, A_kA_l > a_ka_l, k,l \in \{1,2,3\}$.

Then $T$ can be isometrically $PL$ embedded in
$\mathbb{R}^{3}$, as the pleated surface included in the right prism with base
$t$, such that $A_kA_lA_p$ fits $a_kE'_pa_lE'_ka_pE'_l$, where: $B'b  \bot t, B'a_p = R$ and $E'_p, E'_k, E'_l$ on the faces of the
prism, such that $a_kE'_p = Ei_pa_l = \frac{1}{2}A_kA_l$.
%

The following variations of the basic construction above are also considered:
               \begin{enumerate}
               \item Each angle $\varphi$ of $T$ satisfies the condition
$0 < \alpha < \varphi$ and $C\cdot A_kA_l > a_ka_l\,, C < 1$. Moreover,
$A_kA_l/a_ka_l \approx 1$.
               \item Each of the lateral faces of the prism -- including the broken lines
$a_kE'_pa_l$ -- can be (independently) slightly rotated around the lines
$a_ka_l$ such that the construction still can be performed. (The rotation angle depends upon the
constants $\alpha$ and $C$ above.)
               \end{enumerate}

               (In general, one has to simultaneously construct a large number of the standard construction elements above.)
       \end{enumerate}

This employment of the basic construction element may be done straightforwardly in the region of ``almost conformality'' and in the neighborhoods of vertices $v_i$ for which $\theta_i < 2\pi$.
    For vertices $v_i$, such that $\theta_i > 2\pi$, a certain variation of the construction is needed: more precisely $PF$ ``ripples'' are added, see \cite{bz} (and \cite{Sa10} for a short presentation). Note that here, the construction of a (local) triangulation by acute triangles is essential. 

\item We obtain a $\mathcal{C}^2$ smooth, short\footnote{again, in precise, controlled manner} embedding $f_2(\mathcal{P})$ of all the standard construction elements, therefore achieving the desired $PL$ isometric embedding of $\mathcal{P}$.


\end{enumerate}


\section{Obstructions to the Burago-Zalgaller construction in dimension $n \geq 3$}

We list here a number of impediments in the extension of the Burago-Zalgaller construction to higher directions.

\begin{enumerate}

\item {\it Existence of acute triangulations.}

As we have underlined in the previous section (e.g. in (6), (7)), the proof of \cite{bz}, is based on the existence of acute triangulations, more specifically, on a previous result of Burago and Zallgaler \cite{bz0}.\footnote{As we have see, they are necessary in the definition of the canonical construction element, in the construction modified neighborhoods of the vertices, as well as in the triangulation of the compliment of the neighborhoods of the vertices (section 3, (3)). In consequence, they are needed in the final modification of the construction (Section 3, (6)).}
However, such a result does not exists for higher dimensions. Indeed,
next to nothing is known about the existence of
such triangulations in dimension $n \geq 3$.\footnote{There seems to exist little information apart from the one summarized in 
 \cite{za}. (However, there is a renewed hope, due to a different
method recently developed by Tasmuratov \cite{Ta}.)}

\item {\it Nonexistence of smooth embeddings.}

As mentioned in (2) of the previous section, a smooth -- i.e. of class $C^2$ or higher -- embedding of the given $PL$ manifold is supposed to exist. However, the existence of obstructions for the smoothening of a $PL$ manifolds are classical (see \cite{mun1}, \cite{hm}). It follows that, in certain cases, even the first step of the Burago-Zalgaller proof can not be implemented. 

\item {\it Nonexistence of $PL$ approximations.}

As it is shown in 
the proof's synopsis 
above (see also the enunciation of Theorem 1.4 of \cite{bz}), the gist of the proof is to produce $PL$ isometric embeddings arbitrarily close to a given smooth one. In fact, the approximation is even {\it piecewise flat}\footnote{This fact makes the Burago-Zalgaller construction apparently intuitive and attractive for the Graphics community (see \cite{Sa10}).}, at least away from Kuiper waves (cf. Section 3, (4)).

However, such approximations (neither $PL$, nor piecewise flat) do not always exist in codimension 2, as it is shown in a number of counterexamples due to Shtan'ko \cite{Sh}. (It should be noted that they do exist, however, in any other codimension -- see \cite{Lu}.)

Since we work in codimension 1, this obstruction is not truly relevant. It is, however, an impediment if one tries to apply the original Nash construction \cite{na1}, that makes appeal to codimension 2, that is without using Kuiper's improvement \cite{ku}, for which only one additional dimension is needed.

\item {\it Standard conformal mapping}

In Section 3, (2) the role of the standard conformal mapping is described. Furthermore, the mapping $f_2$, that is supposed to be almost conformal 
is introduced in  Section 3, (4). While the authors do not explicitly state this fact, they introduce a {\it quasiconformal structure} on the given manifold, that is further deformed to become arbitrarily close to {\it conformality}.
%

Unfortunately, while any topological manifold of any dimension $n \neq 4$ admits a quasiconformal structure,  by a result of Donaldson and Sullivan \cite{DS}, this does not hold in dimension $n = 4$.
%
%
In fact, there exists an embedding of the unit ball $\mathbb{B}^4$ into $\mathbb{R}^4$, that admits no quasiconformal approximation. It follows that, in dimension $n=4$, the use of the standard conformal mapping is problematic.
Moreover, while $PL$ and locally flat quasiconformal approximations of embeddings exist for $n \geq 2, n \neq 4$ and codimension $m \neq 2$, (see Luukkainen \cite{Lu}), such approximations do not exists for $n = 4$ and $q = 2$ (again, by Shtan'ko's results).


\item {\it Main obstruction}

    However, the main obstruction in extending the Burago-Zalgaller construction to dimension 3 and higher, resides in the fact, already alluded to 
    in Remark \ref{rem:bound-dist}, that it is not possible to obtain the needed ``almost conformal'' quasiconformal parametrization of the neighborhoods of the vertices of a polyhedron of dimension $n \geq 3$, due to the fact that the dilatation is bounded away from $1$ as  discussed in some detail in Section \ref{sec:coefCQ} above. This fact represents the main tool employed in the proof of our main results, in Section 5 below.

 \begin{rem}
 It should be noted that even the shape of the initial smooth embedding is important. Indeed, general ``apple-shaped'' domains in $\mathbb{R}^n, n \geq 3$, do not admit any quasicomformal mapping onto $\mathbb{B}^n$ (see \cite{car}, \cite{vuo1}).
 \end{rem}

\end{enumerate}



\section{Proofs of the Main Results}


Before proceeding to the technical part of the proofs, let us note that, as their name suggests, quasiconformal mappings represent, indeed, the proper, technical notion for the ``almost conformal'' mappings of \cite{bz} -- see \cite{Ag} for precise estimates as well as proofs of the more delicate aspects of the theory, regarding the non-differentiability of everywhere of quasiconformal mappings. (Similar estimates were also given in \cite{Pes}.)

\subsection*{Proof of Theorem \ref{thm:main}}
We show that if $n \geq 3$ one can not produce, for any given $PL$ manifold, quasiconformal embeddings(mappings, as close to conformality as desired,  as it is required in the Burago-Zalgaler construction.


First, let us note that the mapping $f_k$ obtained after performing each of the iterations of the Burago-Zalgaller construction is, indeed, quasiconformal: By Section 3, (4) it is quasiconformal in the complement of the neighborhoods of the vertices. Moreover, in the neighborhoods of the vertices, the construction of Section 3, (3) renders, by the finiteness of the triangulation and by the finite types of simplices employed, a quasiconformal\footnote{even if not, as already stressed, actually conformal} mapping. Therefore, the resulting mapping is, piecewise quasiconformal\footnote{in a rather strong sense}. It follows, by \cite{v0}, that it is, in fact, quasiconformal.

    Moreover, since the polyhedron is compact, the (Nash-Kuiper) limiting isometric $f_2$ mapping will also be quasiconformal (see, e.g. \cite{v0}, Theorem 37.2) 
    and, by  \cite{MRV}, Lemma 4.5, its index will be $\geq \limsup_{k \rightarrow \infty}i(v_j,f_k)$, where $v_j$ denote vertices 
    of the $PL$ manifold.

Let $F_1,F_2$ be two $n$-dimensional faces, $F_1 \cap F_2 = e$, and let $T_e \subset \mathcal{P}$ be a tube of radius $\epsilon$.
Let $\alpha = \measuredangle(F_1,F_2)$ denote the dihedral angle between $F_1,F_2$. The dihedral wedge $D_\alpha$
is raylike at any interior point of $e$ and $T_e \subset D_\alpha$ satisfies the conditions of Theorem
\ref{thm:ray}. It follows that $K(T_e) \geq K(D_\alpha) = \frac{\pi}{\alpha} \gg 0$. Therefore, $K(\mathcal{P})
\geq \max\frac{\pi}{\alpha}$, where the maximum is taken over all the dihedral angles of $\mathcal{P}$.

Clearly, one can produce dilatation $K(\mathcal{P})$ as large as
desired, by choosing polyhedra with at least one dihedral angle
(between $n$-faces) $\pi/m$, where $m$ is any (arbitrarily large)
natural number.\footnote{Note that, in any case, by Corollary 2.18, the distortion coefficient of any polyhedron of dimension $n \geq 3$ is bounded away from 1.}



\hspace*{11.5cm}{$\Box$}

\begin{rem}
By a 
result of Heinonen and Hinkkanen \cite{HH}, at each stage of the NKBZ construction, the resulting map will not only be quasiconformal, but actually {\it quasisymmetric}\footnote{Recall that, given two metric spaces $(X,d)$ and $(Y, \rho)$, an embedding (in particular, a homeomorphism) $f:X \rightarrow Y$ is called {\it quasisymmetric} iff there exists a homeomorphism $\eta:[0,\infty) \rightarrow [0,\infty)$, such that if $d(x,y) \leq td(x,z)$, then $\rho(f(x),f(y)) \leq \eta(t)\rho(f(x),f(z))$, for any triple points $x,y,z \in X$. Note that, while any quasiymmetric homeomorphism is quasiconformal, the opposite implication is far from being trivial even for mappings from $\mathbb{R}^n$ to itself and holds only if $n \geq 2$. In fact, it is false for $n = 1$. (See, e.g. \cite{v}, \cite{He} for proofs and further details and \cite{He1} for a brief, yet  lucid and inspiring exposition.)}.
That is, for compact polyhedra in some $\mathbb{R}^m$, the local, infinitesimal condition of quasiconformality implies (in fact, it is is actually equivalent to) a global one (i.e. quasisymmetry), thence in this case the local distortion of the $PL$ Riemannian metric translates into a global one, {\it quantitatively}\footnote{That is the numerical parameters obtained depend solely on the parameters presumed in the hypothesis.}. (This fact should be viewed in the light of the brief discussion in Section 2.1).
\end{rem}

\begin{rem}
Note that in the theorem of Rickman and Srebro the quasiregular mapping considered is defined on a proper domain $D$ in $R^n$ (for some $n \geq 2$). As such, it applies to $PL$-solid polyhedra $\tilde{\mathcal{P}}$ in $\mathbb{R}^n$, that is to the interior of a compact polyhedral surfaces $\mathcal{P}$ embedded ($PL$ isometrically) in Euclidean space. This corresponds perhaps to our intuition, but falls somewhat short of our goal. To remedy this deficiency, one possibility is to use the fact that $PL$ quasiconformal mappings are quasiconformal \cite{v0} and classical extension results (see, e.g. \cite{v}, 17.18 and 35.3),\footnote{See also another extension result due to Tukia and V\"{a}is\"{a}l\"{a} \cite{TV}.}
to show that the mapping can be extended from $\tilde{\mathcal{P}}$ to $P = \partial\tilde{\mathcal{P}}$.

\end{rem}

\subsection*{Proof of Proposition 1.4}
Proceeding along lines similar to those of the proof of Theorem \ref{thm:main}, suffices in this case to consider a polyhedron (even, for simplicity, one homotopic to a sphere) whose set of vertices satisfies the density conditions required in the statement of Theorem 2.10, and whose (solid) angles at the said vertices are large enough to ensure that the index of the winding mapping at these these points will satisfy condition (2) of Theorem 2.10.

\hspace*{11.5cm}{$\Box$}

Before concluding this section we bring the following remark:

Given the fact that there exists no other method of $PL$ isometric embedding, apart from  
the {\it NKBZ} construction,
it is possible that there exist no such embedding for $PL$ manifolds, if $n \geq 3$. Therefore, in the light of the results proven here and of the more general discussion in the preceding section, 
we venture 
the following

\begin{conj}
For any $n \geq 3$, there exists an $n$-dimensional $PL$-manifold $M^n$ (and in fact, an infinity of such manifolds), that admits no $PL$ isometric embedding in $\mathbb{R}^{n+1}$. 
\end{conj}


\section{The Role 
of Curvature}

\subsection{First Remarks}



Note that the main obstruction in obtaining an almost conformal
mapping resides -- rather against geometric intuition -- on the
edges of triangulation, and not at the vertices ($0$-dimensional
faces). (See also \cite{pe}, pp. 175 and 186.)
%
%
In particular, if $n=3$ and $M^n$ is a a manifold with boundary
embedded in $\mathbb{R}^n$, then small values of the dihedral angle $\alpha$ (see Section 3 and 5 above), are
associated to large mean curvature $H$ of the $PL$ surface
$S^2 = \partial M^3$, as opposed to the Gauss curvature
concentrated at the vertices (see, e.g. \cite{Ba}, \cite{cms}).
%
%
Indeed, any pyramid with large base angles (i.e. with
corresponding large dihedral angles) can be quasiconformally
mapped onto a half-space, with bounded dilatation, which depends
only on $n$ and not on the angles at the vertices $A_i$, even if
the   (dihedral) face-angles, incident to the apex are not
bounded from below (i.e having small, positive Gauss
curvature) (see \cite{car}, Theorem 3.6.10. and Theorem 3.6.13.).\footnote{Moreover, cones -- and even cylinders, as a limiting case when the vertex angle tends to 0 -- can be quasiconformally
mapped onto $\mathbb{H}^n_+$ with
small dilatation, and this can be done independently of the vertex
angle.}

In dimension $n \geq 3$, mean and Gauss curvature are replaced by the so called {$j^{th}$ {\it mean curvatures} and {\it Lipschitz-Killing curvatures}, respectively (see \cite{cms}). Fittingly, these curvatures are also expressible in terms of (higher dimensional) dihedral angles. However, in order not to diverge too much, for details we refer the reader to the above mentioned paper of Cheeger et al.



\begin{rem}
Actually, one can dispense with the use of the initial smooth embedding, and formulate the bounded curvature condition for the $PL$ ($PF$ surface) in terms of the so called {\it generalized principal curvatures}\footnote{We recall the definition of generalized principal curvatures

First, we define be the {\it reach} of $X \subset \mathbb{R}^n$, ${\rm reach}(X)$ as:
\[{\rm reach}(X) = \sup\{r>0\,|\, {\rm for}\; {\rm all}\; y \in N_r(X), {\rm there}\; {\rm exists}\; {\rm a}\; {\rm unique}\; x \in X, {\rm nearest}\; {\rm to}\; y \}\,,\] %
where $N_r(X)$ denotes the $r${\it -neighborhood} of $X$.

Then, for any $0 < \varepsilon < {\rm reach}(X)$, $\partial N_\varepsilon(X)$ is a $C^{1,1}$ hypersurface, hence it admits principal curvatures (in the classical sense) a.e., and let $\kappa_i^\varepsilon(x + \varepsilon\mathbf{n})$, denote the $i$-th such curvature at a generic point. (Here $\mathbf{n}$ represents the normal to $X$ at $x$.)
We can now define the {\it $i$-th generalized principal curvature} as:

\[\kappa_i(x,\mathbf{n}) = \lim_{\varepsilon \rightarrow 0}{\kappa_i^\varepsilon(x + \varepsilon\mathbf{n})}\,,\]
(and the limit exists $\mathcal{H}^{n-1}$ a.a. $(x,\mathbf{n})$, where, as usual, $\mathcal{H}$ denotes the Hausdorff measure).
}
-- see \cite{Zah1}, \cite{Zah2}.
(The fact that the $\mathcal{C}^2$ smooth requirement is, in fact, too strong is also noted in \cite{bz}, Remark 1.10. Note that a sufficient condition mentioned there is that the initial embedding $f_0$ admits a $\mathcal{C}^0$ approximation; compare with the discussion in Sect 4,\,(3).)
\end{rem}

\subsection{Unbounded 2-Manifolds}

We remark here that the NKBZ construction actually holds even for unbounded (noncompact), orientable 2-manifolds, that admit a ``geometrically well behaved'' isometric embedding in $\mathbb{R}^3$. Formally, we can state the following result:

\begin{prop}
Let $M$ be a connected, oriented 2-dimensional $PL$ manifold,
without boundary or having a finite number of compact boundary components, that admits a $\mathcal{C}^2$ smooth embedding into $R^3$.
Then $M$ admits a  $PL$ isometric embedding in $R^3$.
\end{prop}
%

\begin{rem}
Regarding the existence of an isometric embedding, with the prescribed curvature property see Section 6.3 below.
\end{rem}

\begin{proof}
We begin with the simpler case of manifolds without boundary.

Note that, again, by Nash's argument (see \cite{na1}) we may presume that the given $\mathcal{C}^2$ embedding is short.

We can apply Peltonen's argument \cite{pel} (after eventually considering a smoothing to the $\mathcal{C}^\infty$ class) to produce an exhaustion
$\{M_i\}$ of $M$.\footnote{Incidentally, Peltonen's goal was to produce a {\it thick} (or {\it fat}) triangulation of each of the elements of the exhaustion, as well as of their intersections. While more technical definitions can be given (see \cite{Sa05}), for our present purposes suffices to say that thick triangulations are precisely those for which the individual simplices may each be mapped onto a
standard (Euclidean) $n$-simplex, by a $L$-bilipschitz map, followed by a homotety, with a
fixed $L$.}
The ``size'' of the elements of the exhaustion\footnote{as well as the (local) density of the vertices of the triangulation } is a function solely of the {\it maximal osculatory} ({\it tubular}) {\it radius} $\omega_M = \sup\{\rho > 0\,|\, \mathbb{S}^{2}(x,\rho) \; {\rm osculatory}\; {\rm at\; any}\; x \in M\}$, where $\mathbb{S}^{2}(x,\rho)$ is an {\rm osculatory sphere} at $x \in M$ iff: (i) $\mathbb{S}^{\nu-1}(x,\rho)$ is tangent at x; and (ii) $\mathbb{B}^n(x,\rho) \cap M^n = \emptyset$. (For details see \cite{pel}.) Note that an osculatory sphere exists at any point of $M^n$, for all sufficiently smooth Riemannian manifolds -- see, e.g. \cite{pel}. We exploit this feature  to ensure that
that also in this case, the manifold does not intersect itself, hence that an embedding of a tubular neighborhood of $M$ can be obtained, thus assuring that canonical construction elements can be produced without intersections.
%
%

The principal curvatures will be uniformly bounded on each $M_i$, thence, by Buragogo and Zalgaller's argument (\cite{bz}, p. 379)), the $NKBZ$  construction can be 
performed on any of these submanifolds (with boundary) $M_i$.
Therefore, the respective maps $f_i$ and the resulting map $f$
will be quasiconformal, by \cite{Ge}. Here we have to make appeal to the variation for compact manifolds with boundary of the Burago-Zalgaller construction.\footnote{Note that the gluing condition of the pieces $M_i$ is also given Peltonen in terms of $\omega_M$.}
Also, we may have to use a subdivisions and $\varepsilon$-{\it moves}\footnote{see, e.g., \cite{mun}} to assure that, for all $i$, the vertices on the common boundary $N_i = \partial M_i = \partial M_{i+1}$ appertain to the triangulations of both of the considered elements of the exhaustion.

To obtain the extension of the proposition to manifolds with boundary, one has to use the modification of Peltonen's construction given in \cite{Sa05} and apply the Burago-Zalgaller construction for each of the (compact) boundary components.

\end{proof}

\begin{rem}
The fitting result for nonorientable 2-manifolds can also be proved along the same lines.
\end{rem}

\subsection{Metric Curvature and Compatibility Conditions}

We begin by noting that the frustrating, and not infrequently confusing aspect of $PF$ embeddings (and, in general, of $PL$ ones), is that they are highly counterintuitive, not least with regard to the discrepancy between the {\it intrinsic} and {\it extrinsic} curvature. We have discussed this in some detail in \cite{Sa10}. However, we mention here an example due to O'Rourke \cite{O'R}, of a vertex in piecewise flat surface, for which the  intrinsic Gaussian curvature, computed using the {\it angle defect} at the vertex (see, e.g. \cite{Ba}) is $0$, while its extrinsic one, given either by its generalized principal curvatures (as above) or, alternatively, using metric curvatures (see below) is highly positive.

Since the classical notion of curvature can not be defined for $PL$ ($PF$) surfaces, due to their lack of differentiability (at the vertex points), an analogue of Nash's Theorem for smooth 2-manifolds is not immediate, and even the use of generalized curvatures (as mentioned above) does not really solve the problem, and not solely due to the aforementioned counterintuitiveness.

We propose here an approach that allows us to formulate both local and global ``curvature sensitive'' embedding conditions in $\mathbb{R}^3$ for $PL$ 2-manifolds.
Before proceeding further, we should emphasize here that the method we propose here is quite different from the one of, e.g., \cite{Wu}, \cite{BS}, \cite{No}, \cite{MTW}.
We make appeal to notions and results from {\it metric geometry}. In this we rest mainly 
on \cite{Pl}.
First, we introduce some 
notation:

Let $S^n_\kappa$ denote the $n$-dimensional simply connected {\it space form} (i.e. $S^n_0 \equiv \mathbb{R}^n; S^n_\kappa \equiv \mathbb{S}^n_{\sqrt{\kappa}}$ -- the $n$-dimensional sphere of radius $\sqrt{\kappa}$, if $\kappa > 0$; and $S^n_\kappa \equiv \mathbb{S}^n_{\sqrt{\kappa}}$ -- the $n$-dimensional $S^n_{\kappa} \equiv
\mathbb{H}^n_{\sqrt{-\kappa}}$ stands for the hyperbolic space of
curvature $\sqrt{-\kappa}$, as represented by the Poincar\'{e}
ball model of radius $R = 1/\sqrt{-\kappa}$\,, if $\kappa < 0$).

Given a metric space $(X,d)$ and $x,y,z$ points in $X$. The {\it triple} $\{x,y,z\}$ (viewed as a finite metric space) can be isometrically embedded in $S^n_\kappa$ if
\begin{equation}
d(x,y) + d(x,z) + d(y,z) \leq 2\pi\,.
\end{equation}
where $\sqrt{-\kappa}$ is taken to be $\infty$, if $\kappa \leq 0$.

The image of such an isometric embedding is unique, up to an isometry of $S^n_\kappa$ (see, e.g. \cite{Pl}, Proposition 12), it is called the {\it model} (or {\it representative}) {\it triangle} (of the triple $\{x,y,z\}$) in $S^n_\kappa$, and will be denoted by $T(\bar{x},\bar{y},\bar{z})$.

Given three points $x_i,x_j,x_l$ in a metric space $(X,d)$, we denote by $\alpha_\kappa(x_i,x_j,x_l) \in [0,\pi]$, the angle $\measuredangle(x_jx_ix_l)$\footnote{that is of apex $x_i$} of the model triangle in $S^n_\kappa$.
Note that $\alpha_\kappa$ is a continuous function of $\kappa$ and a monotone increasing function in the variables $\kappa$ and $d_1 = d(x_j,x_l)$.

Let $Q = \{x_1,x_2,x_3,x_4\}$ be a {\it metric quadruple} (that is, a 4 point metric space $Q$). 
We introduce the following 
quantities:
\begin{equation}
V_\kappa(x_i) = \alpha_\kappa(x_i;x_j,x_l) + \alpha_\kappa(x_i;x_j,x_m) + \alpha_\kappa(x_i;x_l,x_m)\,
\end{equation}
where $x_i,x_j,x_l,x_m \in Q$ are distinct , and $\kappa$ is any number;
\begin{equation}
A_\kappa(Q) = \max_i{V_\kappa(x_i)}\,.
\end{equation}

From the analogous properties of $\alpha_\kappa$ it follows that $V_\kappa(x)$ is a continuous, monotone increasing function of $\kappa$. (Hence $A_\kappa(Q)$ is also monotone increasing as a function of $\kappa$.)

The first result we need to introduce our approach 
to a curvature-sensitive isometric embedding of $PL$ manifolds is the following proposition (see, e.g. \cite{Pl}, Proposition 20):

\begin{prop}
Let $Q = \{x_1,x_2,x_3,x_4\}$ be a nondegenerate\footnote{A metric quadruple is called {\it nondegenerate} iff no point lies between two other points.} metric quadruple. Then $Q$ admits an isometric embedding in $S^3_\kappa$ iff (i) $A_k(Q) \leq 2\pi$, and (ii) the triangle inequality holds for any triple $\alpha_\kappa(x_i;x_j,x_l)$, $\alpha_\kappa(x_i;x_j,x_m)$, $\alpha_\kappa(x_i;x_j,x_p)$.
Moreover, the embedding is {\em planar} (i.e. it($Q$) can be embedded in (some) $S^2_\kappa$) if there exists an index $i$ such that $\alpha_\kappa(x_i;x_j,x_l) = \alpha_\kappa(x_i;x_j,x_m) + \alpha_\kappa(x_i;x_j,x_p)$.
\end{prop}

The notion on which our approach rests is the so called {\it Wald-Berestovskii} ({\it metric}) {\it curvature} \cite{Wa}, \cite{Ber} (see also \cite{Bl} for a detailed exposition on metric curvatures in general):

\begin{defn} \label{def:WBcurv}
Let $(X,d)$ be a metric space. 
An open set $U \subset X$ is called a {\it region of curvature} $\geq \kappa$ iff any metric quadruple can be isometrically embedded in $S_m$, for some $m \geq k$.
A metric space $(X,d)$ is said to have(to be of) {\it Wald-Berestovskii curvature} $\geq \kappa$ iff for any $x \in X$ is contained in a region $U$ of curvature $\geq \kappa$.
\end{defn}

It turns out that regions of curvature $\geq \kappa$ can be characterized easily in terms of the embedding angle, more precisely we have the following result (see, e.g. \cite{Pl}, Theorem 23):

\begin{prop}
Let $(X,d)$ be a metric space and let $U \in X$ be an open set. $U$ is a region of curvature $\geq \kappa$ iff $V_\kappa(x) \leq 2\pi$, for any metric quadruple $\{x,y,z,t\} \subset U$.
\end{prop}

Note that we can consider the Wald-Berestovskii curvature at an accumulation point (of a metric space) by considering limits of the curvatures of (nondegenerate) regions of diameter converging to 0.
Since, by a theorem of Wald \cite{Wa}, for smooth surfaces (in $\mathbb{R}^3$), Gauss curvature and Wald-Berestovskii coincide,
we can now proceed and present 
our approach to the isometric embedding of $PL$ manifolds into $\mathbb{R}^{n+1}$ problem. First, the metric space we shall consider will be the 1-skeleton of the manifold, with the obvious metric. We note that, in this context the natural choice for the open set $U$ required in Definition \ref{def:WBcurv} is the closed {\it star} of a given vertex $v$, that is, the set $\{e_{vj}\}_j$ of edges incident to 
$v$. Therefore, the set of metric quadruples 
containing the vertex $v$ is finite and Proposition 6.7 is readily applied. 


The local isometric embedding condition is now easy to express, in view of 
Propositions 6.5 and 6.7. Namely, given a vertex $v$, the following system of inequalities:
\begin{equation} \label{eq:metricGauss}
 \left\{
         \begin{array}{lll}
         \max{A_0(v)} \leq 2\pi;\\
         \alpha_0(v;v_j,v_l) \leq \alpha_0(v;v_j,v_p) + \alpha_0(v;v_l,v_p), & {\rm for\; all\;} v_j,v_l,v_p \sim v;\\ 
         V_\kappa(v) \leq 2\pi;
         \end{array}
 \right.
\end{equation}
(Here ``$\sim$'' denotes incidence, i.e. the existence of a 
connecting edge $e_{i} = vv_j$, and, of course, $V_\kappa(v) = \alpha_\kappa(v;v_j,v_l) + \alpha_\kappa(v;v_j,v_p) + \alpha_\kappa(v;v_l,v_p)$, where $v_j,v_l,v_p \sim v$, etc.)

Here the first two inequalities represent the (extrinsic) embedding condition, while the third one represents the intrinsic curvature (of the $PL$ manifold) at the vertex $v$.

The global embedding condition follows immediately: 

\begin{equation} \label{eq:metricNash}
 \left\{
         \begin{array}{lll}
         \max{A_0(v_i)} \leq 2\pi;\\
         \alpha_0(v_i;v_j,v_l) \leq \alpha_0(v;v_j,v_p) + \alpha_0(v;v_l,v_p), & {\rm for\; all\;} v_j,v_l,v_p \sim v_i; \\ 
         V_\kappa(v_i) \leq 2\pi; {\rm for\; all\;} v_i \in V_M.
         \end{array}
 \right.
\end{equation}
where the inequalities above hold for all $v_i \in V_M$. (Here $V_M$ denotes the set of vertices of $M$.)

The system (\ref{eq:metricGauss}) represents 
as a $PL$ 
version of the {\it Gauss compatibility} (or {\it fundamental}) {\it equation} of classical Differential Geometry of surfaces, while (\ref{eq:metricNash}) functions as a $PL$ analogue of the similar global conditions on curvature that are satisfied in the classical Nash embedding of smooth manifolds.
However, for the problem of the path isometric embedding itself, perhaps other methods (such as those adopted in the papers mentioned at the begining of this subsection) should be better considered.
Moreover, this should be regarded in view of the recent results of Matou\v{s}ek, Tancer, and Wagner, \cite{MTW}, namely that the $PL$ isometric embedding problem of $n$-dimensional simplicial complexes in $R^n$, $n \geq 5$ is undecidable and also that the more general problem of $PL$ isometric embedding problem of $n$-dimensional simplicial complexes in $R^N$, is $NP$-hard for $N \geq n \geq (2N - 2)/3$, for any $n \geq 4$.


At this point, a number of 
observations 
are mandatory:

\begin{enumerate}

\item Obviously there exists an inherent weakness in the approach above, due to the fact that the Wald-Berestovskii curvature is a {\it comparison curvature}, thence only inequations, 
    can be given, not equations.
    However, in the defense of the considered definition (of curvature), it may be said that it requires only simple computations, using just  quite standard, elementary trigonometry (albeit in $S^n_\kappa$, for $n = 2$ and $n = 3$).

\item  This brings up the following natural question: Is it possible -- and if so, how? - to actually compute the Wald-Berestovskii curvature of a $PL$ manifold (or of a metric graph), using 
    solely the metric of the manifold (respectively, graph), that is without making recourse to an actual embedding?

The answer to this question is positive, at least for 
spaces satisfying the local existence of shortest geodesic  and having bounded (pinched) curvature \cite{Ber}, Theorem 6 (see also \cite{Bl})
 -- thence also in our case -- it is possible to compute the embedding curvature of a metric quadruple, due to a pioneering work of Wald \cite{Wa} (see also \cite{Bl} for a somewhat more recent and detailed exposition, and \cite{Sa06} for a perhaps more readily available ``digest''). We have the following formula for the embedding curvature $\kappa(Q)$ of a metric quadruple $Q$:

 \begin{equation} \label{eq:k(Q)}
 \kappa(Q) = \left\{
         \begin{array}{clclcrcr}
           \mbox{0} &  \mbox{if $D(Q) = 0$\,;} \\
           \mbox{$\kappa,\, \kappa < 0$} & \mbox{if $det({\cosh{\sqrt{-\kappa}\cdot d_{ij}}}) = 0$\,;} \\
           \mbox{$\kappa,\, \kappa > 0$} & \mbox{if $det(\cos{\sqrt{\kappa}\cdot d_{ij}})$ and $\sqrt{\kappa}\cdot d_{ij} \leq
           \pi$}\\
           & \mbox{\,\, and all the principal minors of order $3$ are $\geq 0$;}
         \end{array}
   \right.
\end{equation}
where $d_{ij} = d(x_i,x_j), 1 \leq i,j \leq 4$, and $D(Q)$ denotes the so called {\it Cayley-Menger determinant}:

\begin{equation}                                \label{eq:D}
 D(x_1,x_2,x_3,x_4) = \left| \begin{array}{ccccc}
                                            0 & 1 & 1 & 1 & 1 \\
                                            1 & 0 & d_{12}^{2} & d_{13}^{2} & d_{14}^{2} \\
                                            1 & d_{12}^{2} & 0 & d_{23}^{2} & d_{24}^{2} \\
                                            1 & d_{13}^{2} & d_{23}^{2} & 0 & d_{34}^{2} \\
                                            1 & d_{14}^{2} & d_{24}^{2} & d_{34}^{2} & 0
                                      \end{array}
                               \right|\;.
\end{equation}

As far as the 
actual computation of $\kappa(Q)$ using Formula (\ref{eq:k(Q)}) is concerned, it should be noted that, apart from the Euclidean case, the equations involved are transcendental, and can not be solved, in general, using elementary methods. Moreover, when solving them by with the assistance of computer based methods (e.g. making use of {\tt MATLAB}), they display certain numerical instability. For a more detailed discussion and some first numerical results, see \cite{Sa04}, \cite{sa}. (Here, again, the advantage of the approach suggested by Propositions 6.5 and 6.7 is evident, at least as far as the type of the involved computations is concerned.)

\item  To compute the extrinsic curvature of an embedding (and to compare it to the intrinsic, given one) we can again make appeal to metric curvatures. In this case, we shall use metric versions for the curvatures of curves, to compute the maximal and minimal sectional curvatures of the embedding. Again, the metric space under investigation is the 1-skeleton of the $PL$ surface, and the considered curves are pairs of edges having in common the vertex at which we wish to compute curvature. As options for the metric curvature of such $PL$ curves we can consider either the {\it Menger curvature} or the {\it Finsler-Haantjes curvature} -- see \cite{Sa04} and \cite{sa} for a discussion of the the practicability of this approach in Graphics, and for some first experimental results (as well as \cite{Bl} for a detailed presentation of the two types of curvatures considered).

\item  While due to the monotony of $V_k$, it would appear that for higher $\kappa$ there are less possible solutions for the third inequality in (\ref{eq:metricNash}) and (\ref{eq:metricGauss}), it should be remembered that this inequality represents a prescribed condition, representing the curvature of the manifold at the considered vertex.

\item To actually solve 
 the system (\ref{eq:metricNash}) appears to be quite difficult and we postpone this problem for further study. We should however, note here that here we consider the isometric embedding (under curvature constrains) of $PL$ manifolds, and not the $PL$ isometric embedding of such manifolds, as in Burago and Zalgaller paper (as well as in the first part of this paper). Paradoxically, the introduction of ``superfluous'' vertices, and, in consequence, 
 of additional inequalities, may actually make the problem more manageable, by introducing more degrees of freedom, thus rendering the embedding more flexible. We do not, however, know how to quantize this freedom and introduce the 
 equations.

\item 

    Obviously, there is nothing special about the particular (and rather restricted) case of 
    1-skeleta of polyhedra (or, more specific, of $PL$ surfaces), and the embedding criteria (\ref{eq:metricGauss}) and (\ref{eq:metricNash}) can be applied to any locally finite metric graph.
    Therefore, system (\ref{eq:metricNash}) represents a partial answer to a question of posed by N. Linial \cite{Li} (see \cite{Sa10}).
    Admittedly, 
    in this case the notion of intrinsic curvature is less natural. Moreover, in the absence of triangles -- whose existence is not guaranteed in the general -- the very definition of, say, Finsler-Haantjes curvature is, 
     in this case, quite problematic. We defer such 
     problems for further study \cite{Sa-prep}.

\item It is sometimes desirable (especially in an applicative setting) to embed the given $PL$ surface (or graph) either in $\mathbb{S}^3$ (for representation reasons), or in $\mathbb{H}^3$ (due to its exponential volume growth). In both cases, the appropriate embedding conditions are immediate to obtain, from Proposition 6.5(i).

    Let us also note here that an isometric embedding criterion in $\mathbb{R}^n$, as well as in $\mathbb{S}^n$ and $\mathbb{H}^n$, in terms of the Cayley-Menger and related determinants (and some related results) can be found in \cite{Bl}. Here, again, the practicability of the required numerical computations is far from clear, and the advantage of our method is evident. 



\end{enumerate}


%
%


\section*{Acknowledgments}

The author wishes to express his gratitude to 
Peter Buser for his wonderful
hospitality and for the numerous stimulating discussions on this
theme and many others that took place during the author's 
visits at EPFL.
The author also would like to thank Konrad Polthier for making him aware 
the paper of Burago and Zalgaller, Nati Linial for his interest and a number of inspiring discussions, and Joseph O'Rourke, for the example mentioned 
in Section 6.3 and for bringing to his attention the related paper of Zalgaller \cite{za2}.
Thanks are also due to David Gu, for his keen interest and warm encouragement.



\end{document}